\documentclass[12pt]{amsart}
\usepackage{geometry}                % See geometry.pdf to learn the layout options. There are lots.
\usepackage{graphicx}
\usepackage{amssymb}
\usepackage{epstopdf}
\usepackage{amsthm}
\newtheorem{theorem}{Theorem}
\newtheorem{lemma}{Lemma}
\newtheorem{corollary}{Corollary}
\newtheorem{proposition}{Proposition}
\newtheorem{definition}{Definition}
\newcommand{\holder}[1]{\|#1\|_\frac12}

\title{A Simple Search Problem}
\author{Marshall W. Buck}
\address{Center for Communications Research}
\author{Douglas H. Wiedemann}

\begin{document}

\begin{abstract}
A simple problem is studied in which there are $N$ boxes and a prize known to be in one of the boxes. Furthermore, the probability that the prize is in any box is given. It is desired to find the prize with minimal expected work, where it takes one unit of work to open a box and look inside. This paper establishes bounds on the minimal work in terms of the $p=1/2$ H\"older norm of the probability density and in terms of the entropy of the probability density. We also introduce the notion of a ``Cartesian product'' of problems, and 
determine the asymptotic behavior of the minimal work for the $n$th power of a problem. 

(This article is a newly typeset version of an internal publication written in 1984. The second author passed away on November 12, 2020, and his estate has approved the submission of this paper.)
\end{abstract}
\maketitle

%\section{}
%\subsection{}
\section{A Simple Problem}
 
Suppose a prize is located in one of $N$ boxes, labeled 1 through $N$. You are told for each box $i$ the probability $p(i)$ that the prize is in that box. We refer to $p(i)$ as the \emph{hiding density}.   What is the strategy which minimizes the expected number of boxes to be opened before obtaining the prize? Since one can remember the boxes that have already been opened the best strategy will be to arrange the numbers from 1 to $N$ into some sequence $a_1, a_2, \ldots ,a_N$, opening box $a_1$ first. If the prize is not immediately found, box $a_2$ is opened, and so on. The expected number of boxes opened (the expected work) before finding the prize is 
\begin{equation}
W = \sum_{j=1}^N j p(a_j).
\end{equation}
Here we demand that the work of opening the box containing the prize must be done even if it is known to be in that box. It is obvious that the best strategy is to open the most likely box first, the next most likely second, and continue in this manner. If there are two boxes which have the same probability, then it 
will not matter which is opened first. Similar search problems are studied and analogs to the following Theorem~\ref{thm:1} are proved in \cite{Kadane80,Kadane78}.

\begin{theorem}\label{thm:1}
The strategy which achieves the minimal number of expected box openings is to order the boxes from $a_1$ to $a_N$ so that $p(a_1) \geq p( a_2) \geq  \cdots \geq p( a_N) $ and then open the boxes in this order until the prize is found. 
\end{theorem}

This minimal value of the work (the ideal work) we denote by $W_I$. The 
ideal strategy is only ideal in terms of number of boxes opened. Implementing it apparently requires sorting the probabilities. The sorting work may be prohibitive when $N$ is very large. Even if we are willing to do the sorting work it is not 
easy to estimate what the value of $W_I$ will be before doing the sort. 

\section{A Simple Strategy}

Consider the following memory-less strategy, which sometimes opens boxes that have already been opened. Fix some probability density $q(1), q(2), \ldots,q(N)$ on the boxes. Select a box $i$ randomly, according to the \emph{search density} $q$, then open box $i$. Continue making random selections until the prize is found. This we will call the \emph{random} strategy. Although this strategy may seem ill-conceived, it will provide bounds on $W_I$. 
First, it is necessary to find the best density $q$. If the prize is actually in box $i$ then the expected number of box openings until the prize is found is $1/q(i)$ 
because this is the expected waiting time for an event which has probability $q(i)$ at each step; i.e., 
$$
1q(i)+2(1-q(i))q(i)+3(1-q(i))^2q(i)+ \dots = \frac1{q(i)}.
$$ 
Thus, we find that the expected work under this strategy is 
\begin{equation}\label{eqn:2}
W= \sum_{i=1}^N \frac{p(i)}{q(i)}.
\end{equation}

It may be assumed that each $p(i)>0$, for there is no need to include boxes with zero probability in the problem. Considering the $(N-1)$-dimensional simplex of possible $q$ densities, it is clear that a minimum of $W$ is achieved in the interior because $W$ goes to infinity as we approach the boundary of the simplex. Applying Lagrange multipliers to locate interior extrema of (2) under the constraint that the sum of $q$ is 1, gives the following equation: 
$$
\frac{p(i)}{q(i)^2}+ \lambda = 0. 
$$
Note the $q(i)$ variables must be proportional to the square roots of the $p(i)$ variables, normalized to make $q$ a probability distribution. The interior minimum is therefore unique. The minimum work for this random strategy can now be found by substitution into (2) and we have the following result. 

\begin{theorem}\label{thm:2}
The minimum work achievable by the random strategy is attained only when 
$$
q(i) = \frac{\sqrt{p(i)}}{\sum_j \sqrt{p(j)}}.
$$
This minimum work is 
$$
\left( \sum_{i=1}^N \sqrt{ p(i)}\right)^2 = \| p \|_{\frac12}.
$$
\end{theorem}

\begin{corollary}\label{cor:to2}
$W_I \leq \| p \|_{\frac12}$.
\end{corollary}

Our next goal will be to examine the relationship between the ideal work and the best random work in more detail. 

\section{How Good Is The Random Strategy?}\label{sec:3}

It will now be assumed that the prior probabilities have been sorted so that $p(1)\geq p(2)\geq \cdots\geq p(N)$. Then the ideal work is $W_I = \sum ip(i)$. A useful technique in comparing this to the expected work of the best random strategy is to invent other random strategies and then use the fact that they have expected work at least as large as $\| p \|_{\frac12}$, the value for the best random strategy.

For example, consider a random strategy with $q(i)$ proportional to $i^{-1}$.  The normalizing factor is $H_N = \sum_{i=1}^N i^{-1}$, the $N$th harmonic number. The expected  
work under this random strategy is by \eqref{eqn:2},  
$$
H_N\cdot \sum_{i=1}^N \frac{p(i)}{i^{-1}} = H_N\cdot W_I.
$$

\begin{theorem}\label{thm:3}
\end{theorem}
$$
\frac{ \| p \|_{\frac12} }{1+\log N} \leq \frac{ \| p \|_{\frac12} }{H_N}\leq W_I \leq  \| p \|_{\frac12}.
$$

\begin{proof}
A random strategy was demonstrated which has work $W_I\cdot H_N$.  Since this is at least as much work as the best random strategy we have 
$\| p \|_{\frac12} \leq W_I\cdot H_N$, 
which gives us the middle inequality in the statement of the theorem. The first inequality now follows from the fact that 
$$
H_N = 1 + \sum_{i=2}^N \frac1i \leq 1 + \int_1^N \frac{dx}{x} = 1 + \log N.
$$
This proves the two leftmost inequalities in the statement of the theorem. The rightmost part is Corollary~\ref{cor:to2}.  
\end{proof}

The lower bound $H_N^{-1}$ on $W_I/\| p \|_{\frac12}$ is sharp. To show this, for a given $N$ select $p(i)$ proportional to $i^{-2}$, $1<i <N$. Then the random search with $q(i)$ proportional to $i^{-1}$ is in fact the best random search and the ratio of ideal work to the work of the best random strategy is precisely $H_N^{-1}$. 

In any random strategy there is the potential for duplicated work. For any strategy, let $D$ be the expected number of times we open a box that has already been opened once, before finding the prize. Were there some way of marking the boxes as we open them, the best strategy would only involve $W-D$ expected work. Note that $W-D$ is an upper bound on the ideal work. 
\begin{theorem}\label{thm:4}
For the best random strategy $D = (W-1)/2$. 
\end{theorem}
\begin{proof}
Suppose the prize is actually in box $i$. For $j\neq i$ we can compute the expected number of times box $j$ is ``wastefully" opened before opening box $i$ and finding the prize, where a wasteful opening of a box is any opening beyond the first one. If our random strategy is defined by a density $q$ then all that matters is the relative values of $q(i)$ and $q(j)$. For the immediate purpose neglect all box openings except those of box $i$ or $j$. The result is a sequence of Bernoulli trials with probability $r = q(j)/( q(i)+q(j))$ that one opens box $j$. The expected number of wasteful openings of box $j$ is 
$$
1r^2(1-r)+ 2r^3(1-r) + \cdots = \frac{r^2}{1-r} = \frac{q(j)^2}{q(i)(q(i)+q(j))}.
$$
Summing over all pairs $j\neq i$, we obtain

\begin{eqnarray*}
D &=& \sum_i p(i) \sum_{j\neq i} \frac{q(j)^2}{q(i)(q(i)+q(j))} =
\sum_{i,j,i\neq j} p(i) \frac{q(j)^2}{q(i)(q(i)+q(j))}\\
  &=& \sum_{i,j} p(i) \frac{q(j)^2}{q(i)(q(i)+q(j))}- \frac12.
\end{eqnarray*}
This last step is justified because 
$$
\sum_{i} p(i) \frac{q(i)^2}{q(i)(q(i)+q(i))}= \frac12.
$$

So far in the proof we have not used the specific values of $q(i)$. In fact, for the best strategy we know $q(i)$ is proportional to $\sqrt{p(i)}$.  Since the function $q(j)^2/(q(i)( q(i)+q(j)))$ is homogeneous of degree zero the normalizing factor 
doesn't matter. For the best random strategy we have 
\begin{eqnarray*}
\sum_{i,j} p(i) \frac{q(j)^2}{q(i)(q(i)+q(j))} &=& 
\sum_{i,j} p(i) \frac{p(j)}{\sqrt{p(i)}(\sqrt{p(i)}+\sqrt{p(j)})}\\
&=& \sum_{i,j} \frac{p(j)\sqrt{p(i)}}{\sqrt{p(i)}+\sqrt{p(j)}}.
\end{eqnarray*}
The last expression can be symmetrized with respect to $i$ and $j$ by adding the sum to itself with $i$ and $j$ reversed and dividing by two. The result is 
\begin{eqnarray*}
\frac12 \sum_{i,j} \frac{p(j)\sqrt{p(i)}+p(i)\sqrt{p(j)}}{\sqrt{p(i)}+\sqrt{p(j)}}
&=&\frac12 \sum_{i,j} \sqrt{p(i)}\sqrt{p(j)}\\
&=& \frac12 \left( \sum_i \sqrt{p(i)} \right)^2 = \frac12 \| p \|_\frac12 = \frac{W}{2}.
\end{eqnarray*}
Thus, we have shown $D = (W-1)/2$. 
\end{proof}

\begin{corollary}\label{cor:to4}
$$ W_I \leq (\| p \|_\frac12+1)/2.
$$
\end{corollary} 
\begin{proof}
The bound in Theorem~\ref{thm:2} can be improved by subtracting the duplicate work (Theorem~\ref{thm:4}) of the best random strategy.  
\end{proof}

\section{The Connection With Entropy}\label{sec:4}

The previous section established bounds on the ideal work in terms of $\| p \|_\frac12$ where $p$ is the hiding density. The ideal work is obviously a measure of the ``spread" of the density $p$. A traditional measure of the spread of $p$ is the entropy: 
$$
h(p) = -\sum_{i=1}^N  p(i) \log p(i). 
$$
Actually the entropy is roughly the \emph{logarithm} of the number of likely positions, so we will compare the ideal work to $e^{h(p)}$.
First, it will be shown that $e^{h(p)}$ is never more than the work of the best random strategy. 
\begin{theorem}\label{thm:5}
$e^{h(p)}\leq \| p \|_\frac12$.
\end{theorem}
\begin{proof}
First note that we have: 
$$
\log \| p \|_\frac12 = 2\log \sum_i p(i)^\frac12 = 2 \log \sum_i p(i) p(i)^{-\frac12}.
$$
Then, by the convexity of the logarithm function, we obtain: 
$$
2 \log \sum_i p(i) p(i)^{-\frac12}\geq 2  \sum_i p(i) \log (p(i)^{-\frac12}) = h(p).
$$
\end{proof}

Now we find the hiding distributions on $N$ points with a given ideal work $W_I$ and maximum entropy. These will turn out to be densities $p$ where $p(i)$ is proportional to $x^i$ for some $x>0$. We therefore let 
$$
p_x(i) = p_{x,N}(i) = \frac{x^i}{\sum_{j=1}^N x^j}.
$$
Furthermore, let $M_x = \sum_{i=1}^N i p_x(i)$. 

\begin{lemma}\label{lemma:1} Let $N> 1$ be a positive integer and let $W$ be any real number such that $1<W<N$. There is a unique $x>0$ such that $M_x= W$. 
\end{lemma}
\begin{proof} 
The first thing to prove is that for $N> 1$, $M_{x,N}$ is a strictly increasing function of $x$. We show this by induction. For $N=2$, we have $M_x = 2 -1/(1+x)$, which is a strictly increasing function. Assume that it is true for $N$ and we are to prove it for $N+ 1$. For $N> 1$ we have 
$$
M_{x,N+1} = M_{x,N} + \phi(x)(N+1-M_{x,N}) = (1-\phi(x))M_{x,N}+ \phi(x)(N+1)
$$
where
$$
\phi(x) = \frac{x^{N+1}}{\sum_{j=1}^{N+1} x^j} = \left( \sum_{j=0}^N x^{-j}\right)^{-1} < 1.
$$
Then
$$
\frac{d}{dx} M_{x,N+1} = (1-\phi) \frac{d}{dx}M_{x,N} + (N+1-M_{x,N})\frac{d\phi}{dx}.
$$
The first summand is positive by the inductive hypothesis. The second summand is positive because $\phi(x)$ is an increasing function of $x$ and because $M_{x,N}< N+1$, which in turn is true because $M_{x,N}$ is the mean value of a density on integers less
than $N+1$.
 
The lemma follows because $M_{x,N}$ is a continuous function which approaches 1 as $x\to 0$ and approaches $N$ as $x\to \infty$.
\end{proof}
 
That $M_{x,N}$ is an increasing function of $x$ is also a consequence of the more general result below (proof omitted).
\begin{proposition}\label{prop:1}
Suppose that a probability measure $\mu$ has compact support on the reals $\mathbf{R}$. Then the function defined by 
$$
F(z) = \int e^{xz} d\mu(x)
$$
is real-analytic and log-convex. 
\end{proposition}

Taking $\mu$ above to be $\sum_{k=1}^N \delta_k$, we get $F(z) = \sum_{k=1}^N e^{kz}$. Now, $F$ is log-convex if and only if the function $F'/F$ is increasing. However, the latter function is just $M_{e^x,N}$, so we also know that $M_{x,N}$ is increasing. 

\begin{proposition}\label{prop:2}
Let $N> 1$ be a positive integer and let $W$ be any real number such that $1< W <N$. Then, of those densities on the set $\{1,2,\ldots, N\}$ which have mean $W$ there is a unique one of maximum entropy and it is the density of the form $p_{x,N}$ which has mean W.
\end{proposition}
\begin{proof} 
By Lemma~\ref{lemma:1} there is an $a>0$ such that $p_a$ has mean $W$.  Let $q$ be any density with mean $W$. Our proof makes use of the Kullback-Leibler inequality, 
$$
\sum_{i=1}^N q(i) \log \frac{q(i)}{p_a(i)} \geq 0
$$
with equality if and only if $q(i) = p_a(i)$ for all $i$ in the range $1 \leq i \leq N$. Thus, 
$$
\sum_{i=1}^N q(i) \log q(i)+\sum_{i=1}^N q(i)\log\frac{\sum_{i=1}^N a^i}{a^i}\geq 0
$$
and we must have: 
\begin{equation}\label{eqn:3}
h(q) \leq -W\log a+ \log \sum_{i=1}^N a^i.
\end{equation}
Equality holds if and only if $q=p_a$; we have for $q \neq p_a$, $ h(q) < h(p_a)$. 

Note, for future reference, that \eqref{eqn:3} holds for all $a>0$, not only for that one value for which $p_{a,N}$ has mean $W$. 
\end{proof}

Now we seek a \emph{lower} bound on the ideal work $W$ in terms of the entropy of the hiding density. (Later we show that it is impossible to find a reasonable upper bound on the ideal work in terms of the entropy.) 

\begin{theorem}\label{thm:6}
If $p$ is a hiding density on a finite number of boxes and its ideal work is $W$, then $h(p) \leq \log( W^W/(W-1)^{(W-1)})$ and $e^{h(p)-1} < W$. 
\end{theorem}
\begin{proof}
Starting as we did in the proof of Proposition~\ref{prop:2}, we obtain, for 
\emph{every} $a>0$: 
$$
h(p) \leq -W\log a+ \log \sum_{i=1}^N a^i.
$$
Restricting ourselves to $a<1$, we can replace the right-hand side by the infinite sum and still have a valid inequality. 
$$
h(p) \leq (W-1) \log \left(\frac1a\right) + \log \left(\frac1{1-a}\right).
$$
Assume that $W>1 $. The function on the right-hand side is continuous on the 
open interval $(0,1)$, approaches $+\infty$ on each end, and achieves its minimum value when $a = (W-1)/ W$. Substituting this choice for $a$ and simplifying, we obtain 
$$
h(p) \leq \log( W^W/(W-1)^{W-1}),
$$
which is the first inequality in the statement of the Theorem. 

The second inequality follows from the first. Apply the Mean Value Theorem to the function $f(x) = x \log x$ and the two endpoints $W-1$ and $W$: there exists an $x_0$ such that $W-1 < x_0 < W$ and such that 
$$
\frac{df}{dx} (x_0) = \frac{f(W)-f(W-1)}{W-(W-1)}.
$$
This expands to 
$$
1 + \log x_0 = W\log W - (W-1)\log (W-1) = \log( W^W/(W-1)^{W-1}). 
$$
But $1+\log x_0 < 1+\log W$, so by combining inequalities we obtain $h(p) < 1+\log W$, and consequently $e^{h(p)-1} < W$. 

When $W=1$, there can only be one box to open, so $h(p)=0$ , and the two inequalities still hold. 
\end{proof}

Theorem~\ref{thm:6} cannot be improved to read $e^{h(p)-b} \leq W$ for a value of $b$ less than 1. Counterexamples to any such assertion can be obtained by using the distributions $p_{x,\infty}$, letting $x$ approach $1$ from below. In fact, fixing $x<1$ and letting
$N$ grow, we see that $W_{x,\infty} = \lim_{N\to\infty} W_{x,N} = 1/(1-x)$. Furthermore, we can compute
the limiting entropy: 
$$
h(p_{x,\infty}) = \frac{x}{1-x}\log \frac1x + \log\frac1{1-x}.
$$
Then, 
$$
\lim_{x\to 1^-} \{ h(p_{x,\infty})- \log W_{x,\infty} \} =  
\lim_{x\to 1^-} \frac{x}{1-x}\log \frac1x  = 1.
$$

It is natural to ask if there are upper bounds on the ideal work in terms of the entropy alone. There is no such bound as can be seen by considering the hiding densities on $N$ points with $p(i)$ made proportional to $i^{-2}$, studied in Section~\ref{sec:3}. The entropies of these hiding densities remain bounded as $N\to \infty$ because 
$$
-\sum_{i=1}^\infty \log(i^{-2}) i^{-2}
$$
has a finite sum. 
The ideal work, however, goes to infinity with $N$ because 
$$
\sum_{i=1}^\infty i\cdot i^{-2}
$$
diverges. 

\section{Cartesian Products}\label{sec:5}
 
Given that $p$ and $q$ are densities on finite sets $A$ and $B$, let $p\times q$ be the ``Cartesian product'' density on $A\times B$. This is defined by $p \times q((a, b)) = p(a)q(b)$ for each $(a, b)\in A \times B$. It is not clear how the ideal work behaves under this product but the behavior of our bounds follows from the well known: 
\begin{eqnarray}
\label{eqn:a} \holder{p\times q} &=& \holder{p} \holder{q}\\
\label{eqn:b}h(p\times q)&=& h(p)+h(q).
\end{eqnarray}

The purpose of this section is to study the ideal work of the Cartesian product of a density with itself many times, $p \times p \times \cdots \times p = p^{\times k}$. From Theorem~\ref{thm:3} we can obtain the basic exponential growth rate: 
$$
\lim_{k\to\infty} \frac{\log W(p^{\times k})}{k} = \log \holder{p}.
$$
However, in order to obtain an asymptotic result of the form $$W(p^{\times k}) = (1+o(1))g(k)\holder{p}^k,$$ we must work harder. Our main device will be a map $\phi$ which takes a density on a set of $N$ objects to a density on at most $N$ real numbers. Specifically, 
\begin{definition}
For a density $p$ on a finite set $A$ let 
$$
\phi(p) = \sum_{a\in A} p(a) \delta_{\log p(a)}
$$
where $\delta_x$, for $x\in \mathbf{R}$, denotes the probability measure on $\mathbf{R}$ which gives probability 1 to the point $x$.
\end{definition}
 
Note $\phi(p)$ is written as a measure on $\mathbf{R}$, but it can also be viewed as a probability density on $\{ \log p(a)\mid a \in A \}$. In addition it will be convenient to define a map $\psi$ which is like $\phi$, except that the delta functions have constant coefficients. 
\begin{definition}
If $p$ is a density on a finite set $A$ of cardinality $N$ let 
$$
\psi(p) = \sum_{a\in A} N^{-1} \delta_{\log p(a)}.
$$
\end{definition}

In order to relate these definitions to the ideal work we need the following result. 
\begin{lemma}\label{lemma:2}
For any hiding density $p$ on the finite set $\{1,2,\ldots,N\}$,
$$
W_I(p) = \sum_{\substack{i,j\\p(i)<p(j)}} p(i) + \frac12\sum_{\substack{i,j\\p(i)=p(j)}} p(i) + \frac12.
$$
\end{lemma}

\begin{proof}
It may be assumed without loss of generality that 
$p(1)\geq p(2)\geq \cdots \geq p(N)$. Then, 
\begin{eqnarray*}
W_I = \sum_{i=1}^N i p(i) &=& 
\sum_{\substack{i,j\\p(i)<p(j)}} p(i) + \sum_{\substack{i,j \\ j<i \\ p(i)=p(j)}} p(i) +\sum_i p(i)\\
&=& \sum_{\substack{i,j\\p(i)<p(j)}} p(i) + \frac12 \sum_{\substack{i,j \\ i\neq j \\ p(i)=p(j)}} p(i) +\sum_i p(i)\\
&=& \sum_{\substack{i,j\\p(i)<p(j)}} p(i) + \frac12 \sum_{\substack{i,j \\ p(i)=p(j)}} p(i) +\frac12 \sum_i p(i)\\
&=& \sum_{\substack{i,j\\p(i)<p(j)}} p(i) + \frac12 \sum_{\substack{i,j \\ p(i)=p(j)}} p(i) +\frac12.
\end{eqnarray*}
\end{proof}

The relevance of these definitions to the ideal work is made clear in the following theorem. Given any measure $\alpha$ on $\mathbf{R}$ let $\alpha^R$ denote the measure reflected 
about the origin. Also let the convolution of two probability measures 
$\alpha,\beta$ on $\mathbf{R}$ be denoted by $\alpha*\beta$.

\begin{theorem}\label{thm:7}
Define the measure  $\mu = \phi(p)^R* \psi(p)$. Then
$$
W_I = N\mu(\{x\mid x>0\})+ \frac{N}2\mu(\{0\})+ \frac12.
$$
\end{theorem}
\begin{proof} 
We can again assume that $p$ is a density on $\{1,2,\ldots,N \}$. Corresponding terms in the statement of the Theorem~\ref{thm:7} and the Lemma~\ref{lemma:2} are equal. Consider the first term: 
\begin{eqnarray*}
\phi(p)^R* N\psi(p)(\{x\mid x>0\})&=& \sum_{i,j} p(i) \cdot \delta_{\log p(j) - \log p(i)}(\{x\mid x>0\})\\
&=& \sum_{\substack{i,j\\p(i)<p(j)}} p(i).
\end{eqnarray*}
The second terms are equal by a similar argument.
\end{proof}

Recall that the goal of this section is to investigate $W_I(p^{\times k})$ for large $k$. The following properties are all easy consequences of the definitions just made. 
\begin{eqnarray}
\phi(p\times q) &=& \phi(p)* \phi(q) \label{eqn:c}\\
\psi(p\times q) &=& \psi(p)* \psi(q) \label{eqn:d}\\
(\alpha*\beta)^R &=& \alpha^R*\beta^R \label{eqn:e}
\end{eqnarray}

In light of Theorem~\ref{thm:7} we seek to find 
$(\phi(p)^R* \psi(p))^{* k}$ where $* k$
denotes convolution of the measure with itself $k$ times. How can the measure of the positive half line be approximated? The central limit theorem shows that in a sense this can be approximated by a Gaussian density with mean and variance that are $k$ 
times the mean and variance of $\mu = \phi(p)^R* \psi(p)$. Unfortunately, this Gaussian approximation will in general give a poor estimate of $\mu^{* k}(\{x\mid x>0\})$ because the mean of $\mu$ may be far from zero. 

A technique which saves the day is ``Chernoff tilting''\cite{Chernoff}.     If a probability measure $\alpha$ of bounded support is multiplied pointwise by the function $e^{sx}$, another measure results. This measure can be normalized, resulting in a new (``tilted'') probability measure. The idea of Chernoff tilting is that the new measure is a version of $\alpha$, but has its mean shifted to a new place. Tilting $\mu$ so that the resulting measure has zero mean will turn out to help a great deal. First there are a few technicalities. 

\begin{definition}
If $\alpha$ is a probability measure of compact support on $\mathbf{R}$ and $s$ is any real number, let $e^{sx}\alpha$ be the measure defined by $e^{sx}\alpha(A) = \int_A e^{sx} d\alpha(x)$. Define 
the \emph{total} of this measure be $t(\alpha,s) = \int e^{sx} d\alpha(x)$. Finally, define the \emph{tilted} probability measure to be $\tau(\alpha,s) = t(\alpha,s)^{-1} e^{sx}\alpha$. 
\end{definition}

These definitions obey the following properties
\begin{eqnarray}
\label{eqn:f}  e^{sx}(\alpha*\beta) &=& (e^{sx}\alpha)* (e^{sx}\beta) \\
\label{eqn:g}  t(\alpha*\beta,s) &=& t(\alpha,s)\cdot t(\beta,s) \\
\label{eqn:h}  \tau(\alpha,s)^R  &=&  \tau(\alpha^R,-s)\\
\label{eqn:i}  \tau(\tau(\alpha,s),t) &=&  \tau(\alpha,s+t)\\
\label{eqn:j}   \tau(\alpha*\beta,s) &=&  \tau(\alpha,s)* \tau(\beta,s).
\end{eqnarray}

We don't provide proofs of these results but they all follow directly from the 
definitions. What tilt of $\mu = \phi(p)^R* \psi(p)$ should be taken? Note that 
$$
e^x\psi(p) = \sum_{a\in A} N^{-1} e^{\log p(a)} \delta_{\log p(a)} =
\sum_{a\in A} N^{-1} p(a) \delta_{\log p(a)}.
$$
Thus, 
\begin{equation}
\label{eqn:k} \tau(\psi(p),1) = \phi(p).
\end{equation}
\begin{lemma} \label{lemma:3}
If $\mu = \phi(p)^R* \psi(p)$ then
$\tau(\mu,1/2) = \tau(\psi(p),1/2)^R* \tau(\psi(p),1/2)$ and
$t(\mu,1/2) = N^{-1}\holder{p}$.
\end{lemma}
\begin{proof} Use \eqref{eqn:j},\eqref{eqn:h},\eqref{eqn:k},\eqref{eqn:i} to derive the sequence of equalities: 
\begin{eqnarray*}
\tau(\phi^R* \psi,1/2) &=& \tau(\phi^R,1/2)* \tau(\psi,1/2) =
    \tau(\phi,-1/2)^R*\tau(\psi,1/2)\\
    &=& \tau(\tau(\psi,1),-1/2)^R* \tau(\psi,1/2) = \tau(\psi,1/2)^R* \tau(\psi,1/2).
\end{eqnarray*}
Also 
\begin{eqnarray*}
t(\mu,1/2)&=& t(\phi^R* \psi,1/2) = t(\phi^R,1/2)t(\psi,1/2) = t(\phi,-1/2)t(\psi,1/2)\\
&=& \left( \sum p(i)^\frac12\right)\left(\sum N^{-1} p(i)^\frac12\right) = N^{-1} \holder{p}.
\end{eqnarray*}
\end{proof} 
What has just been shown is that $\tau(\mu,1/2)$ is the convolution of a probability density $\tau(\psi,1/2)$ with its reverse and is therefore symmetric about zero. Thus $1/2$ is the correct amount by which to tilt. 

\begin{theorem}\label{thm:8}
Let $\zeta$ denote the measure $\tau(\mu,1/2)=\tau(\psi,1/2)^R* \tau(\psi,1/2)$.  Then the ideal work for the ``Cartesian power'' hiding density is given by 
$$
W_I(p^{\times k}) = \frac{\holder{p}^k}{2}\int e^{-|x|/2} d\zeta^{* k}(x) + \frac12.
$$
\end{theorem}
\begin{proof}
Using $e^{x/2}\mu^{* k} = t(\mu,1/2)^k \zeta^{* k}$, it follows that
$$
  \mu^{* k}(\{0\}) = t(\mu,1/2)^k \zeta^{* k}(\{0\})=
  N^{-k} \holder{p}^k \zeta^{* k}(\{0\})
$$
and 
\begin{eqnarray*}
\mu^{* k}(\{x \mid x>0\})&=&
  t(\mu,1/2)^k (e^{-x/2}\zeta^{* k})(\{x \mid x>0\})\\
  &=& N^{-k} \holder{p}^k \int_{(0,\infty)} e^{-x/2} d\zeta^{* k}(x).
\end{eqnarray*}
Applying Theorem~\ref{thm:7} with $N^k$ and $\mu^{* k}$ in place of $N$ and $\mu$, we obtain 
$$
W_I(p^{\times k}) = \holder{p}^k \left(  \frac12 \zeta^{* k}(\{0\})+
 \int_{(0,\infty)} e^{-x/2} d\zeta^{* k}(x)\right)+ \frac12.
$$
Now, because $\zeta$ is the convolution of a measure with its reverse, it is symmetric about zero. Then, for the expression in large brackets above we can substitute 
$$
\frac12 \int e^{-|x|/2} d\zeta^{* k}(x).
$$
\end{proof} 

A good approximation for the integral in Theorem~\ref{thm:8} can be obtained by replacing the density $\zeta^{*k}$ by a normal approximation. This is the Gaussian with 
mean zero and variance $k\sigma^2$ where $\sigma^2$ is the variance of 
$\zeta$.  The validity of this approximation depends on the accuracy of the normal approximation. It appears that if we wish to obtain an asymptotically correct formula for the work then something stronger than the central limit theorem is required. If the measure $\zeta$ is supported in a lattice then we need to use a ``local limit theorem''

\begin{definition}
A measure $\zeta$ on $\mathbf{R}$ is said to be of \emph{lattice type} if there are two real numbers $a$ and $b$ such that $\zeta$ is supported on the set $\{a+bm \mid m\in \mathbf{Z}\}$. 
\end{definition}

If $\zeta$ is supported on a discrete lattice then so are all $\zeta^{*k}$. This happens, for example, when $p$ gives only two boxes nonzero probability. In the case of the symmetric measure $\zeta$ occurring in Theorem~\ref{thm:8}, the set $\{0\}$ has positive measure, so if $\zeta$ were of lattice type the constant $a$ in the definition could be set to 0. In this case, assume that $p$ is not a uniform distribution, so that $\zeta$ is not entirely supported at 0. Let $b$ assume its largest possible positive value. In this case, by the local limit theorem (for the normal distribution) on page 233 of Gnedenko and Kolmogorov\cite{GK}, we approximate the integral in Theorem~\ref{thm:8} asymptotically by the summation
$$
\sum_{m=-\infty}^\infty e^{-|mb|/2} \frac{b}{\sqrt{2\pi k}\sigma} e^{-m^2b^2/2\sigma^2 k}
$$ 
with error of order $o(k^{-1/2})$. The last exponential factor is approximately unity for large $k$. 
Now, 
$$
\sum_{m=-\infty}^\infty e^{-|mb|/2} \frac{b}{\sqrt{2\pi k}\sigma}=
\frac{b}{\sqrt{2\pi k}\sigma\tanh(b/4)}.
$$

In the non-lattice case we use the theorem of Cram\'er and Ess\'een, appearing as Theorem 2, page 210, in \cite{GK}. Let $\nu$ be a non-lattice measure having mean 0, variance $\sigma^2$, and finite third moment $\alpha_3$. Let $\Phi$ denote the normal distribution function: 
$$
\Phi(x) = \frac1{\sqrt{2\pi}} \int_{-\infty}^x e^{-y^2/2} dy.
$$
Define the normalized distribution function 
$$
F_k(x) = \int_{(-\infty,x\sigma\sqrt{k})} d\nu^{*k}(y)
$$
Then the theorem of Cram\'er and Ess\'een states that 
$$
F_k(x) - \Phi(x) = R_k(x) = \frac{\alpha_3}{6\sigma^3\sqrt{2\pi k}}(1-x^2)e^{-x^2/2} + o(k^{-1/2}),
$$
uniformly in $x$. Define the function $F(x,y)$ to be $\frac{e^{-y/2}}2$ on the set $\{ (x,y) \mid y>0\text{\,and\,} -y\leq x<y \}$, and to be zero elsewhere. If we take Lebesgue measure in $y$ and the measure $\nu^{*k}$ in $x$ and apply Fubini's Theorem we obtain the equality of 
$$
\int\int F(x,y) dy d\nu^{*k}(x) = \int e^{-|x|/2} d\nu^{*k}(x)
$$
and 
\begin{eqnarray*}
\int\int F(x,y) d\nu^{*k}(x)dy &=& \frac12\int e^{-y/2} \left( F_k\left(\frac{y}{\sigma\sqrt{k}}\right) - F_k\left(-\frac{y}{\sigma\sqrt{k}}\right)\right) dy\\
&=& \frac12\int e^{-y/2} \left( \Phi\left(\frac{y}{\sigma\sqrt{k}}\right) - \Phi\left(-\frac{y}{\sigma\sqrt{k}}\right)\right) dy + o(k^{-1/2}).
\end{eqnarray*}
Meanwhile the right-hand side of the last formula is also obtained if we replace $\nu^{*k}$ by a continuous Gaussian of variance $k\sigma^2$.  This theorem implies that it is asymptotically correct in our case to replace $\zeta^{*k}$ by the Gaussian. The integral thus obtained is the same as the limit as $b\to 0$ of the formula for the lattice case. 
Notice that when $\holder{p} > 1$ the final $1/2$ in Theorem~\ref{thm:8} is unimportant, and this occurs whenever the hiding density gives positive probability to more than one box. These results are summarized in the following theorem. 

\begin{theorem}\label{thm:9}
If the hiding density $p$ takes at least two distinct positive values then as $k\to\infty$  the ideal work $W_I(p^{*k})$ is asymptotic to one of the following: 

\begin{enumerate}
\item Non-lattice case:  $$ \frac{2\holder{p}^k}{\sqrt{2\pi k}\sigma},$$
\item Lattice case: $$ \frac{b\holder{p}^k}{\sqrt{2\pi k}2\sigma\tanh(b/4)},$$
\end{enumerate}
where
$$
\sigma^2 = \text{Variance}(\zeta) = \holder{p}^{-1} \sum_{i,j} \sqrt{p(i)p(j)} \log^2\frac{p(i)}{p(j)}
$$
and in the lattice case $b$ is the period of the lattice generated by the set of numbers $\{\log\frac{p(i)}{p(j)}\}$.
\end{theorem}

Although $\holder{p}$ started its life as a cheap upper bound for the ideal work, in the Cartesian product limit problem it has acquired legitimacy!

\end{document}